\documentclass[10pt]{article}
\begin{document}
\title{ {\bf An Extremal Problem On Potentially $K_{r+1}-H$-graphic
Sequences}
\thanks{  Project Supported by NNSF of China(10271105), NSF of Fujian(Z0511034),
Science and Technology Project of Fujian, Fujian Provincial
Training Foundation for "Bai-Quan-Wan Talents Engineering" ,
Project of Fujian Education Department and Project of Zhangzhou
Teachers College.}}
\author{{Chunhui Lai, Lili Hu}\\
{\small Department of Mathematics, Zhangzhou Teachers College,}
\\{\small Zhangzhou, Fujian 363000,
 P. R. of CHINA.}\\{\small e-mail: zjlaichu@public.zzptt.fj.cn (Chunhui Lai)}}
\date{}
\maketitle
\begin{center}
\begin{minipage}{4.2in}
\vskip 0.1in
\begin{center}{\bf Abstract}\end{center}
 {Let $K_k$, $C_k$, $T_k$, and $P_{k}$ denote a complete graph on $k$
vertices,  a cycle on $k$ vertices, a tree on $k+1$ vertices, and a
path on $k+1$ vertices, respectively. Let $K_{m}-H$ be the graph
obtained from $K_{m}$ by removing the edges set $E(H)$ of the graph
$H$ ($H$ is a subgraph of $K_{m}$). A sequence $S$ is potentially
$K_{m}-H$-graphical if it has a realization containing a $K_{m}-H$
as a subgraph. Let $\sigma(K_{m}-H, n)$ denote the smallest degree
sum such that every $n$-term graphical sequence $S$ with
$\sigma(S)\geq \sigma(K_{m}-H, n)$ is potentially
$K_{m}-H$-graphical.  In this paper, we determine the values of
$\sigma (K_{r+1}-H, n)$ for
    $n\geq 4r+10, r\geq 3, r+1 \geq k \geq 4$ where $H$ is a graph on $k$
    vertices which
    contains a tree on $4$ vertices but
     not contains a cycle on $3$ vertices. We also determine the values of
      $\sigma (K_{r+1}-P_2, n)$ for
    $n\geq 4r+8, r\geq 3$.} \par
\par
\end{minipage}
\end{center}
 \par
 \section{Introduction}
\par

  The set of all non-increasing nonnegative integers sequence $\pi=$
  ($d(v_1 ),$ $d(v_2 ),$ $...,$ $d(v_n )$) is denoted by $NS_n$.
  A sequence
$\pi\epsilon NS_n$ is said to be graphic if it is the degree
sequence of a simple graph $G$ on $n$ vertices, and such a graph $G$
is called a realization of $\pi$. The set of all graphic sequences
in $NS_n$ is denoted by $GS_n$. A graphical sequence $\pi$ is
potentially $H$-graphical if there is a realization of $\pi$
containing $H$ as a subgraph, while $\pi$ is forcibly $H$-graphical
if every realization of $\pi$ contains $H$ as a subgraph. If $\pi$
has a realization in which the $r+1$ vertices of  largest degree
induce a clique, then $\pi$ is said to be potentially
$A_{r+1}$-graphic. Let $\sigma(\pi)=d(v_1 )+d(v_2 )+... +d(v_n ),$
and $[x]$ denote the largest integer less than or equal to $x$. We
denote $G+H$ as the graph with $V(G+H)=V(G)\bigcup V(H)$ and
$E(G+H)=E(G)\bigcup E(H)\bigcup \{xy: x\in V(G) , y \in V(H) \}. $
Let $K_k$, $C_k$, $T_k$, and $P_{k}$ denote a complete graph on $k$
vertices,  a cycle on $k$ vertices, a tree on $k+1$ vertices, and a
path on $k+1$ vertices, respectively. Let $K_{m}-H$ be the graph
obtained from $K_{m}$ by removing the edges set $E(H)$ of the graph
$H$ ($H$ is a subgraph of $K_{m}$).
\par

Given a graph $H$, what is the maximum number of edges of a graph
with $n$ vertices not containing $H$ as a subgraph? This number is
denoted $ex(n,H)$, and is known as the Tur\'{a}n number. This
problem was proposed for $H = C_4$ by Erd\"os [2] in 1938 and
generalized by Tur\'{a}n [15]. In terms of graphic sequences, the
number $2ex(n,H)+2$ is the minimum even integer $l$ such that every
$n$-term graphical sequence $\pi$ with $\sigma (\pi) \geq l $ is
forcibly $H$-graphical. Here we consider the following variant:
determine the minimum even integer $l$ such that every $n$-term
graphical sequence $\pi$ with $\sigma(\pi)\ge l$ is potentially
$H$-graphical. We denote this minimum $l$ by $\sigma(H, n)$.
Erd\"os,\ Jacobson and Lehel [4] showed that $\sigma(K_k, n)\ge
(k-2)(2n-k+1)+2$ and conjectured that equality holds. They proved
that if $\pi$ does not contain zero terms, this conjecture is true
for $k=3,\ n\ge 6$. The conjecture is confirmed in
[5],[10],[11],[12] and [13].
 \par
 Gould,\ Jacobson and
Lehel [5] also proved that  $\sigma(pK_2, n)=(p-1)(2n-2)+2$ for
$p\ge 2$; $\sigma(C_4, n)=2[{{3n-1}\over 2}]$ for $n\ge 4$. Luo [14]
characterized the potentially $C_{k}$ graphic sequence for
$k=3,4,5.$  Lai [7] determined  $\sigma (K_4-e, n)$ for $n\ge 4$.\
Lai [8, 9] determined
    $\sigma (K_{5}-C_{4}, n),\sigma (K_{5}-P_{3}, n)$ and
    $\sigma (K_{5}-P_{4}, n),$ for  $n\geq 5$.
Yin, Li and Mao[17] determined $\sigma(K_{r+1}-e,n)$ for $r\geq 3,$
$r+1\leq n \leq 2r$ and $\sigma(K_5-e,n)$ for $n\geq5$. Yin and
Li[16] gave a good method (Yin-Li method) of determining the values
$\sigma(K_{r+1}-e,n)$ for $r\geq2$ and $n\geq3r^2-r-1$. After
reading[16],  using Yin-Li method
    Yin[18] determined the values $\sigma(K_{r+1}-K_3,n)$ for
    $r\geq3,n\geq3r+5$. Determining $\sigma(K_{r+1}-H,n)$, where $H$
    is a tree on 4 vertices is more useful than a cycle on 4
    vertices (for example, $C_4 \not\subset C_i$, but $P_3 \subset C_i$ for $i\geq 5$).
    So, after reading[16] and [18],  using Yin-Li method
    we prove the following three
theorems.\par

{\bf  Theorem 1.1.} If $r\geq3$ and $n\geq 4r+8$, then
          $\sigma(K_{r+1}-P_2,n)=(r-1)(2n-r)-2(n-r)+2.$
\par

{\bf  Theorem 1.2.} If $r\geq3$ and $n\geq 4r+10$, then $\sigma
(K_{r+1}-T_{3}, n)= (r-1)(2n-r)-2(n-r).$

\par
{\bf  Theorem 1.3.} If $r\geq3, r+1 \geq k \geq 4$ and $n\geq
4r+10$, then $\sigma (K_{r+1}-H, n)= (r-1)(2n-r)-2(n-r),$ where $H$
is a graph on $k$
    vertices which contains a tree on $4$ vertices but
     not contains a cycle on $3$ vertices.
       \par

     There are a number of graphs on $k$
    vertices which
    containing a tree on $4$ vertices but
     not containing a cycle on $3$ vertices (for example, the cycle
      on $k$ vertices, the tree on $k$ vertices, and the
   complete 2-partite graph on $k$ vertices, etc ).\par
\section{Preparations}\par
  In order to prove our main result,we need the following notations
  and results.\par
  Let $\pi=(d_1,\cdots,d_n)\epsilon NS_n,1\leq k\leq n$. Let \par
    $$ \pi_k^{\prime\prime}=\left\{
    \begin{array}{ll}(d_1-1,\cdots,d_{k-1}-1,d_{k+1}-1,
    \cdots,d_{d_k+1}-1,d_{d_k+2},\cdots,d_n), \\ \mbox{ if $d_k\geq k,$}\\
    (d_1-1,\cdots,d_{d_k}-1,d_{d_k+1},\cdots,d_{k-1},d_{k+1},\cdots,d_n),
     \\ \mbox{if $d_k < k.$} \end{array} \right. $$
  Denote
  $\pi_k^\prime=(d_1^\prime,d_2^\prime,\cdots,d_{n-1}^\prime)$,where
  $d_1^\prime\geq d_2^\prime\geq\cdots\geq d_{n-1}^\prime$ is a
  rearrangement of the $n-1$ terms of $\pi_k^{\prime\prime}$. Then
  $\pi_k^{\prime}$ is called the residual sequence obtained by
  laying off $d_k$ from $\pi$.\par
    {\bf Theorem 2.1[16]} Let $n\geq r+1$ and $\pi=(d_1,d_2,\cdots,d_n)\epsilon
    GS_n$ with $d_{r+1}\geq r$. If $d_i\geq 2r-i$ for
    $i=1,2,\cdots,r-1$, then $\pi$ is potentially $A_{r+1}$-graphic.
    \par
    {\bf Theorem 2.2[16]} Let $n\geq 2r+2$ and $\pi=(d_1,d_2,\cdots,d_n)\epsilon
    GS_n$ with $d_{r+1}\geq r$. If $d_{2r+2}\geq r-1$ , then $\pi$ is
    potentially $A_{r+1}$-graphic.
    \par
    {\bf Theorem 2.3[16]} Let $n\geq r+1$ and $\pi=(d_1,d_2,\cdots,d_n)\epsilon
    GS_n$ with $d_{r+1}\geq r-1$. If $d_i\geq 2r-i$ for
    $i=1,2,\cdots,r-1$, then $\pi$ is potentially $K_{r+1}-e$-graphic.
    \par
    {\bf Theorem 2.4[16]} Let $n\geq 2r+2$ and $\pi=(d_1,d_2,\cdots,d_n)\epsilon
    GS_n$ with $d_{r-1}\geq r$. If $d_{2r+2}\geq r-1$ , then $\pi$ is
    potentially $K_{r+1}-e$
    -graphic.
    \par
    {\bf Theorem 2.5[6]} Let $\pi=(d_1,\cdots,d_n)\epsilon NS_n$ and $1\leq k\leq
    n$. Then $\pi\epsilon GS_n$ if and only if  $\pi_k^\prime\epsilon
    GS_{n-1}$.
    \par
     {\bf Theorem 2.6[3]} Let $\pi=(d_1,\cdots,d_n)\epsilon NS_n$
     with even $\sigma(\pi)$. Then $\pi\epsilon GS_n$ if and only if
     for any $t$,$1\leq t\leq n-1$,
     $$\sum_{i=1}^t d_i\leq t(t-1)+\sum_{j=t+1}^n
     min \{t,d_j \}.$$
     \par
     {\bf Theorem 2.7[5]} If $\pi=(d_1,d_2,\cdots,d_n)$ is a graphic
     sequence with a realization $G$ containing $H$ as a subgraph, then
     there exists a realization $G^\prime$ of $\pi$ containing H as a
     subgraph so that the vertices of $H$ have the largest degrees of
     $\pi$.
     \par
     {\bf Lemma 2.1 [18]} If $\pi=(d_1,d_2,\cdots,d_n)\epsilon NS_n$ is potentially
     $K_{r+1}-e$-graphic, then there is a realization $G$ of $\pi$
     containing $K_{r+1}-e$ with the $r+1$ vertices
     $v_1,\cdots,v_{r+1}$ such that $d_G(v_i)=d_i$ for
     $i=1,2,\cdots,r+1$ and $e=v_rv_{r+1}$.
     \par

     {\bf Lemma 2.2  [18]} If $r\geq3$ and $n\geq r+1$, then
          $\sigma(K_{r+1}-K_3,n)\geq(r-1)(2n-r)-2(n-r)+2.$

    \par

\section{ Proof of Main results.} \par

  {\bf Lamma 3.1} Let $n\geq r+1$ and $\pi=(d_1,d_2,\cdots,d_n)\epsilon GS_n$
with $d_r\geq r-1$ and $d_{r+1}\geq r-2$. If $d_i\geq 2r-i$ for
$i=1,2,\cdots,r-2,$ then $\pi$ is potentially $K_{r+1}-P_2$-graphic.
\par
  {\bf Proof.} We consider the following two cases.
  \par
  Case 1: $d_{r+1}\geq r-1.$
  \par
  Subcase 1.1:  $d_{r-1}\geq r+1$. Then $\pi$ is potentially
  $K_{r+1}-e$-graphic by Theorem 2.3. Hence,$\pi$ is potentially $K_{r+1}-P_2$-graphic.
  \par
   Subcase 1.2: $d_{r-1}=r-1$. Then $d_{r-1}=d_r=d_{r+1}=r-1$.
  \par If
  $d_{r+2}=r-1$, then the residual sequence $\pi_{r+1}^\prime=
  (d_1^\prime,\cdots,d_{n-1}^\prime)$
  obtained by laying off $d_{r+1}=r-1$ from $\pi$ satisfies: (1)
  $d_i^\prime=d_i-1$ for $i=1,2,\cdots,r-2,$ (2)
  $d_1^\prime\geq2(r-1)-1$,$ \cdots,$
  $d_{(r-1)-1}^\prime=d_{r-2}^\prime\geq2(r-1)-(r-2),d_{r-1}^\prime
  =d_r,d_{(r-1)+1}^\prime=d_r^\prime
  =d_{r+2}=r-1$. By Theorem 2.1, $\pi_{r+1}^\prime$ is potentially
  $A_{(r-1)+1}$-graphic. Therefore, $\pi$ is potentially
  $K_{r+1}-P_2$-graphic by $ \{d_1-1,\cdots,d_{r-2}-1,d_r,d_{r+2} \}
  = \{d_1^\prime,\cdots,d_r^\prime \}$
  and Theorem 2.7.
  \par
  If
  $d_{r+2}\leq r-2$, then the residual sequence $\pi_{r+1}^\prime=
  (d_1^\prime,\cdots,d_{n-1}^\prime)$
  obtained by laying off $d_{r+1}=r-1$ from $\pi$ satisfies:
  (1)  $d_i^\prime=d_i-1$ for $i=1,2,\cdots,r-2,$ (2)
  $d_1^\prime\geq2(r-1)-1,\cdots,
  d_{(r-1)-1}^\prime=d_{r-2}^\prime\geq2(r-1)-(r-2),d_{r-1}^\prime
  =d_r,$
   $d_{(r-1)+1}^\prime=d_r^\prime=d_{r-1}-1=r-2$.
   By Theorem 2.3, $\pi_{r+1}^\prime$ is potentially
  $K_{(r-1)+1}-e$-graphic. Therefore, $\pi$ is potentially
  $K_{r+1}-P_2$-graphic by $\{d_1-1,\cdots,d_{r-2}-1,d_r,d_{r-1}-1\}
  = \{d_1^\prime,\cdots,d_r^\prime \}$
  and Lemma 2.1.
  \par
  Subcase 1.3:  $d_{r-1}=r$. Then $d_{r+1}=r$ or $r-1$.

  If $d_{r+1}=r$, then $d_{r-1}=d_r=d_{r+1}=r$. The residual sequence
  $\pi_{r+1}^\prime$ satisfies: (1) $d_i^\prime=d_i-1$ for
  $i=1,2,\cdots,r-2$,(2) $d_1^\prime\geq d_1-1\geq 2(r-1)-1,\cdots,$
  $d_{(r-1)-1}^\prime=d_{r-2}^\prime\geq d_{r-2}-1
  \geq 2(r-1)-(r-2)$ and $d_{(r-1)+1}^\prime=d_r^\prime\geq d_r-1=r-1$. By Theorem
  2.1, $\pi_{r+1}^\prime$ is potentially
  $A_{(r-1)+1}$-graphic. Thus, $\pi$ is potentially
  $K_{r+1}-P_2$-graphic by $\{d_1-1,\cdots,d_{r-2}-1 \}
  \subseteq \{d_1^\prime,\cdots,d_r^\prime \}$
  and Theorem 2.7.
  \par
    If $d_{r+1}=r-1$, then $d_r=r-1$ or $r$.
    \par
    If $d_r=r-1$, then $\pi_{r+1}^\prime$ satisfies: (1) $d_i^\prime=d_i-1$ for
  $i=1,2,\cdots,r-1$,(2) $d_1^\prime\geq d_1-1\geq 2(r-1)-1,\cdots,
  d_{(r-1)-1}^\prime=d_{r-2}^\prime=d_{r-2}-1
  \geq 2(r-1)-(r-2)$ and $d_{(r-1)+1}^\prime=d_r^\prime=d_r=r-1$.
  According to Theorem 2.1, $\pi_{r+1}^\prime$ is potentially
  $A_{(r-1)+1}$-graphic. Therefore, $\pi$ is potentially
  $K_{r+1}-P_2$-graphic by $\{d_1-1,\cdots,d_{r-1}-1,d_r \}=
  \{d_1^\prime,\cdots,d_r^\prime \}$
  and Theorem 2.7.
  \par
    If $d_r=r$, then $\pi_{r+1}^\prime$ satisfies: (1) $d_i^\prime=d_i-1$ for
  $i=1,2,\cdots,r-2$,(2) $d_1^\prime\geq d_1-1\geq2(r-1)-1,
  \cdots,d_{(r-1)-1}^\prime=d_{r-2}^\prime=d_{r-2}-1
  \geq2(r-1)-(r-2)$ and $d_{(r-1)+1}^\prime=d_r^\prime=d_{r-1}-1=r-1$.
  By Theorem 2.1, $\pi_{r+1}^\prime$ is potentially
  $A_{(r-1)+1}$-graphic. Therefore, $\pi$ is potentially
  $K_{r+1}-P_2$-graphic by $\{d_1-1,\cdots,
  d_{r-2}-1,d_r,d_{r-1}-1\}=\{d_1^\prime,\cdots,d_r^\prime \}$
  and Theorem 2.7.
  \par
  Case 2: $d_{r+1}\leq r-2$, that is, $d_{r+1}= r-2$.
\par
If $d_{r-1}<d_{ r-2}$, then $\pi_{r+1}^\prime$ satisfies: (1)
$d_i^\prime=d_i-1$ for
  $i=1,2,\cdots,r-2$,(2) $d_1^\prime=d_1-1\geq 2(r-1)-1,
  \cdots,d_{(r-1)-1}^\prime=d_{r-2}^\prime=d_{r-2}-1
  \geq 2(r-1)-[(r-1)-1]$ and $d_{(r-1)+1}^\prime=d_r^\prime=d_r\geq r-1$.
  By Theorem 2.1, $\pi_{r+1}^\prime$ is potentially
  $A_{(r-1)+1}$-graphic. Therefore, $\pi$ is potentially
  $K_{r+1}-P_2$-graphic by $\{d_1-1,\cdots,d_{r-2}-1,d_{r-1},d_r \}=
  \{d_1^\prime,\cdots,d_r^\prime \}$
  and Theorem 2.7.
  \par
  If $d_{r-1}=d_{ r-2}\geq r+2$, then $\pi_{r+1}^\prime$ satisfies:
   $d_1^\prime \geq d_1-1 \geq 2(r-1)-1,\cdots,d_{(r-1)-1}^\prime=
   d_{r-2}^\prime \geq d_{r-2}-1
  \geq 2(r-1)-[(r-1)-1]$ and $d_{(r-1)+1}^\prime=d_r^\prime \geq r-1$.
  By Theorem 2.1, $\pi_{r+1}^\prime$ is potentially
  $A_{(r-1)+1}$-graphic. Therefore, $\pi$ is potentially
  $K_{r+1}-P_2$-graphic by $\{d_{r-1},d_r,d_1-1, \cdots,d_{r-2}-1 \}=
  \{d_1^\prime,
  \cdots,d_r^\prime \}$
  and Theorem 2.7.
  \par
  {\bf Lemma 3.2.}  Let $n\geq 2r+2$ and $\pi=(d_1,d_2,\cdots,d_n)\epsilon GS_n$
with  $d_{r-2}\geq r$. If $d_{2r+2}\geq r-1$, then $\pi$ is
potentially $K_{r+1}-P_2$-graphic.
\par
  {\bf Proof.} We consider the following two cases.
  \par
  Case 1: If $d_{r-1}\geq r$. Then $\pi$ is potentially
  $K_{r+1}-e$-graphic by Theorem 2.4. Hence, $\pi$ is
potentially $K_{r+1}-P_2$-graphic.
  \par
  Case 2: $d_{r-1}\leq r-1$, that is, $d_{r-1}= r-1$, then
  $d_{r-1}=d_r=d_{r+1}=\cdots=d_{2r+2}=r-1$ and $\pi_{r+1}^\prime$ satisfies: (1)
$d_i^\prime=d_i-1$ for
  $i=1,2,\cdots,r-2$,(2) $d_{(r-1)+1}^\prime=d_r^\prime\geq r-1$ and
   $d_{2(r-1)+2}^\prime=d_{2r}^\prime\geq (r-1)-1$.
  By Theorem 2.2, $\pi_{r+1}^\prime$ is potentially
  $A_r$-graphic. Therefore, $\pi$ is potentially
  $K_{r+1}-P_2$-graphic by $\{d_1-1,\cdots,d_{r-2}-1,d_r,d_{r+2} \}
  =\{d_1^\prime,\cdots,d_r^\prime \}$
  and Theorem 2.7.
  \par
  {\bf Lemma 3.3.} If $r\geq3$ and $n\geq r+1$, then
          $\sigma(K_{r+1}-P_2,n)\geq(r-1)(2n-r)-2(n-r)+2.$
          \par
    {\bf Proof.} By Lemma 2.2, for $r\geq3$ and $n\geq r+1$,
    $\sigma(K_{r+1}-K_3,n)\geq(r-1)(2n-r)-2(n-r)+2.$
    Obviously, for $r\geq3$ and $n\geq r+1$,
      $\sigma(K_{r+1}-P_2,n)\geq
    \sigma(K_{r+1}-K_3,n)\geq(r-1)(2n-r)-2(n-r)+2.$
    \par
    {\bf  Lemma 3.4.} If $r\geq3, r+1 \geq k \geq 4$ and $n\geq
r+1$, then $\sigma (K_{r+1}-H, n)\geq (r-1)(2n-r)-2(n-r),$ for $H$
be a graph on $k$
    vertices which containing a tree on $4$ vertices but
     not containing a cycle on $3$ vertices.
\par
{\bf Proof.}   Let
$$G =K_{r-2}+  \overline{K_{n-r+2}}$$
Then $G$ is a unique realization of $((n-1)^{r-2}, (r-2)^{n-r+2})$
and $G$ clearly does not contain $K_{r+1}-H,$ where the symbol $x^y$
means $x$ repeats $y$ times in the sequence. Thus
$$\sigma (K_{r+1}-H, n)\geq (r-2)(n-1)+(r-2)(n-r+2)+2= (r-1)(2n-r)-2(n-r).$$
\par

    {\bf The Proof of Theorem 1.1 }
       According to Lemma 3.3, it is enough to verify that for
       $r\geq 3$ and
       $n\geq
       4r+8$,
       $$\sigma(K_{r+1}-P_2,n)\leq(r-1)(2n-r)-2(n-r)+2.$$
    We now prove that if $n\geq 4r+8$ and $\pi=(d_1,d_2,\cdots,d_n)\epsilon GS_n$
with
$$\sigma(\pi)\geq(r-1)(2n-r)-2(n-r)+2,$$
then $\pi$ is
potentially
  $K_{r+1}-P_2$-graphic.
  \par
   If $d_{r-2}\leq r-1$, then
   $$ \begin{array}{rcl}
   \sigma(\pi)&\leq &(r-3)(n-1)+(r-1)(n-r+3)\\
   &=&   (r-1)(n-1)-2(n-1)+(r-1)(n-r+3)\\
   &=& (r-1)(2n-r)-2(n-r)\\
   &<& (r-1)(2n-r)-2(n-r)+2,
   \end{array} $$
   which is a
    contradiction. Thus $d_{r-2}\geq r$.
    \par
    If $d_r\leq r-2$, then
    $$ \begin{array}{rcl}
    \sigma(\pi)&=&\sum_{i=1}^{r-1}d_i+\sum_{i=r}^n d_i\\
    & \leq & (r-1)(r-2)+\sum_{i=r}^{n}min \{r-1,d_i \}+\sum_{i=r}^n
    d_i\\
    &=&(r-1)(r-2)+2\sum_{i=r}^n
    d_i\\
    &\leq &(r-1)(r-2)+2(n-r+1)(r-2)\\
    &=& (r-1)(2n-r)-2(n-r)-2\\
    &<& (r-1)(2n-r)-2(n-r)+2,
    \end{array} $$
    which is a
    contradiction. Hence $d_r\geq r-1$.
\par
    If $d_{r+1}\leq r-3$, then
    $$\begin{array}{rcl}
    \sigma(\pi)&=&\sum_{i=1}^{r-1}d_i+d_r+\sum_{i=r+1}^n d_i \\
    & \leq &(r-1)(r-2)+
    \sum_{i=r}^{n}min \{r-1,d_i \}+d_r+\sum_{i=r+1}^n d_i \\
    &=&(r-1)(r-2)+min \{r-1,d_r \}+d_r+2\sum_{i=r+1}^n
    d_i \\
    & \leq & (r-1)(r-2)+2d_r+2\sum_{i=r+1}^nd_i \\

    & \leq &(r-1)(r-2)+2(n-1)+2(n-r)(r-3) \\
    &=&(r-1)(2n-r)-2(n-r)\\
    &<&(r-1)(2n-r)-2(n-r)+2, \end{array} $$
     which is a
    contradiction. Thus $d_{r+1}\geq r-2$.
    \par
      If $d_i\geq 2r-i$ for $i=1,2,\cdots,r-2$ or $d_{2r+2}\geq r-1$, then $\pi$
      is potentially
  $K_{r+1}-P_2$-graphic by Lemma 3.1 or Lemma 3.2. If $d_{2r+2}\leq
  r-2$ and there exists an integer $i$,$1\leq i\leq r-2$ such that
  $d_i\leq2r-i-1$, then
  $$\begin{array}{rcl}
  \sigma(\pi) &\leq &(i-1)(n-1)+(2r+1-i+1)(2r-i-1)\\
  &&+
  (r-2)(n+1-2r-2)\\
  &=& i^2+i(n-4r-2)-(n-1)+(2r-1)(2r+2)\\
  &&+(r-2)(n-2r-1).
  \end{array} $$
  Since $n\geq 4r+8$, it is  easy to see that
  $i^2+i(n-4r-2)$, consider as a function of $i$, attains its maximum
  value when $i=r-2$. Therefore,
  $$\begin{array}{rcl}
  \sigma(\pi) &\leq &
  (r-2)^2+(n-4r-2)(r-2)-(n-1)\\
  &&+(2r-1)(2r+2)+(r-2)(n-2r-1)\\
  &=&(r-1)(2n-r)-2(n-r)+2-n+4r+7 \\
  &<&\sigma(\pi),
  \end{array} $$
  which is a
    contradiction.
    \par
    Thus, $\sigma
(K_{r+1}-P_{2}, n)\leq(r-1)(2n-r)-2(n-r)+2$ for
    $n\geq 4r+8$.
\par

 {\bf The Proof of Theorem 1.2 }
       According to Lemma 3.4, it is enough to verify that for
       $r\geq 3$ and
       $n\geq
       4r+10$,
       $$\sigma(K_{r+1}-T_3,n)\leq(r-1)(2n-r)-2(n-r).$$
    We now prove that if $n\geq 4r+10$ and $\pi=(d_1,d_2,\cdots,d_n)\epsilon GS_n$
with
$$\sigma(\pi)\geq(r-1)(2n-r)-2(n-r),$$
then $\pi$ is potentially
  $K_{r+1}-T_3$-graphic.
  \par
   If $d_{r-2}\leq r-1$, we consider the following cases.
   \par
   (1)Suppose $d_{r-2}= r-1$ and
 $\sigma(\pi)=(r-3)(n-1)+(r-1)(n-r+3)$, then
 $\pi=((n-1)^{r-3},(r-1)^{n-r+3})$. Obviously $\pi$ is potentially
 $K_{r+1}-T_3$ graphic.
 \par
 (2)Suppose $d_{r-2}= r-1$ and
 $\sigma(\pi)<(r-3)(n-1)+(r-1)(n-r+3)$,
 then
 $$ \begin{array}{rcl}
   \sigma(\pi)&<&(r-3)(n-1)+(r-1)(n-r+3)\\
   &=&   (r-1)(n-1)-2(n-1)+(r-1)(n-r+3)\\
   &=& (r-1)(2n-r)-2(n-r),

   \end{array} $$
 which is
 a contradiction.
 \par
 (3)Suppose $d_{r-2}< r-1$, then
   $$ \begin{array}{rcl}
   \sigma(\pi)&< &(r-3)(n-1)+(r-1)(n-r+3)\\
   &=&   (r-1)(n-1)-2(n-1)+(r-1)(n-r+3)\\
   &=& (r-1)(2n-r)-2(n-r),
   \end{array}
    $$
   which is a
    contradiction.
    \par
    Thus, $d_{r-2}\geq r$ or $\pi$ is potentially
 $K_{r+1}-T_3$ graphic.
    \par
    If $d_r\leq r-2$, then
    $$ \begin{array}{rcl}
    \sigma(\pi)&=&\sum_{i=1}^{r-1}d_i+\sum_{i=r}^n d_i\\
    & \leq & (r-1)(r-2)+\sum_{i=r}^{n}min \{r-1,d_i \}+\sum_{i=r}^n
    d_i\\
    &=&(r-1)(r-2)+2\sum_{i=r}^n
    d_i\\
    &\leq &(r-1)(r-2)+2(n-r+1)(r-2)\\
    &=& (r-1)(2n-r)-2(n-r)-2\\
    &<& (r-1)(2n-r)-2(n-r),
    \end{array} $$
    which is a
    contradiction. Hence $d_r\geq r-1$.
\par
    If $d_{r+1}\leq r-3,$ we consider the following cases.
    \par
  (1)Suppose $d_r=n-1$, then $d_1 \geq d_2\geq \cdots \geq d_{r-1}\geq
 d_r=n-1$, therefore $d_1=d_2=\cdots=d_r=n-1$. Therefore $d_{r+1} \geq r,$
 which is a contradiction.
 \par
 (2)Suppose $d_r\leq n-2,$ then
    $$\begin{array}{rcl}
    \sigma(\pi)&=&\sum_{i=1}^{r-1}d_i+d_r+\sum_{i=r+1}^n d_i \\
    & \leq &(r-1)(r-2)+
    \sum_{i=r}^{n}min \{r-1,d_i \}+d_r+\sum_{i=r+1}^n d_i \\
    &=&(r-1)(r-2)+min \{r-1,d_r \}+d_r+2\sum_{i=r+1}^n
    d_i \\
    & \leq & (r-1)(r-2)+2d_r+2\sum_{i=r+1}^nd_i \\

    & \leq &(r-1)(r-2)+2(n-2)+2(n-r)(r-3) \\
    &=&(r-1)(2n-r)-2(n-r)-2 \\
    &<&(r-1)(2n-r)-2(n-r), \end{array} $$
     which is a
    contradiction.
    \par
    Thus $d_{r+1}\geq r-2$.
    \par
      If $d_i\geq 2r-i$ for $i=1,2,\cdots,r-2$ or $d_{2r+2}\geq r-1$, then
      $\pi$ is potentially
 $K_{r+1}-T_3$ graphic($\pi=((n-1)^{r-3},(r-1)^{n-r+3})$) or $\pi$
      is potentially
  $K_{r+1}-P_2$-graphic  by Lemma 3.1 or Lemma 3.2 . Therefore, $\pi$
      is potentially
  $K_{r+1}-T_3$-graphic. If $d_{2r+2}\leq
  r-2$ and there exists an integer $i$, $1\leq i\leq r-2$ such that
  $d_i\leq 2r-i-1$, then
  $$\begin{array}{rcl}
  \sigma(\pi) &\leq &(i-1)(n-1)+(2r+1-i+1)(2r-i-1)\\
  &&+
  (r-2)(n+1-2r-2)\\
  &=& i^2+i(n-4r-2)-(n-1)\\
  &&+(2r-1)(2r+2)+(r-2)(n-2r-1).
  \end{array} $$
  Since $n\geq 4r+10$, it is  easy to see that
  $i^2+i(n-4r-2)$, consider as a function of $i$, attains its maximum
  value when $i=r-2$. Therefore,
  $$\begin{array}{rcl}
  \sigma(\pi) &\leq &
  (r-2)^2+(n-4r-2)(r-2)-(n-1)\\
  &&+(2r-1)(2r+2)+(r-2)(n-2r-1)\\
  &=&(r-1)(2n-r)-2(n-r)-n+4r+9 \\
  &<&\sigma(\pi),
  \end{array} $$
  which is a
    contradiction.
    \par
    Thus, $\sigma
(K_{r+1}-T_{3}, n) \leq(r-1)(2n-r)-2(n-r)$ for
    $n\geq 4r+10$.
\par

 {\bf The Proof of Theorem 1.3 } By Lemma 3.4,
 for $r\geq3, r+1 \geq k \geq 4$ and $n\geq r+1$,
$\sigma (K_{r+1}-H, n)\geq (r-1)(2n-r)-2(n-r).$ Obviously, for
$r\geq3, r+1 \geq k \geq 4$ and $n\geq 4r+10$, $\sigma (K_{r+1}-H,
n)\leq  \sigma(K_{r+1}-T_3, n).$  By theorem 1.2, for $r\geq3, r+1
\geq k \geq 4$ and $n\geq 4r+10$, $\sigma (K_{r+1}-T_{3}, n)=
(r-1)(2n-r)-2(n-r).$  Then $\sigma (K_{r+1}-H, n)=
(r-1)(2n-r)-2(n-r),$  for $r\geq3, r+1 \geq k \geq 4$ and $n\geq
4r+10.$
\par

 \section*{Acknowledgment}
  The authors thanks the referees for many helpful comments.

\par

\end{document}